\documentclass[10pt]{article}

\usepackage[english]{babel}
\usepackage{amsmath}
\usepackage{amssymb}
\usepackage{bm}
\usepackage{bbm}
\usepackage{mathrsfs,color}
\usepackage{graphics,graphicx,theorem}
\usepackage{dsfont}
\usepackage{setspace}

\usepackage{natbib}
\usepackage{multirow}
\usepackage{hhline}
\usepackage{hyperref}

\theoremstyle{plain}
\newtheorem{theorem}{Theorem}[section]

\newtheorem{lemma}{Lemma}[section]

\newtheorem{proposition}{Proposition}[section]

\hypersetup{
	colorlinks=true,
	urlcolor=blue,
	citecolor=blue,
	linktoc=all,
	linkcolor=red} 

\usepackage[scriptsize]{subfigure}
\allowdisplaybreaks

\makeatletter
\renewcommand\@biblabel[1]{}
\makeatother

\def\bSig\mathbf{\Sigma}

\newcommand{\qed}{$\square$}

\usepackage[top=1.2in,bottom=1.2in,left=1.2in,right=1.2in]{geometry}

\begin{document}

\title{\bf {\Large{On consistent estimation of the missing mass}}}
\author{
Fadhel Ayed\\
\normalsize{University of Oxford, UK}\\
\normalsize{email: \texttt{fadhel.ayed@gmail.com}}
\bigskip\\
Marco Battiston\\
\normalsize{University of Oxford, UK}\\
\normalsize{email: \texttt{marco.battiston@stats.ox.ac.uk}}
\bigskip\\
Federico Camerlenghi\\
\normalsize{University of Milano--Bicocca, Italy}\\
\normalsize{email: \texttt{federico.camerlenghi@unimib.it}}
\bigskip\\
Stefano Favaro\\
\normalsize{University of Torino and Collegio Carlo Alberto, Italy}\\
\normalsize{email: \texttt{stefano.favaro@unito.it}}
}
\date{}
\maketitle
\thispagestyle{empty}

\setcounter{page}{1}
\begin{center}
\textbf{Abstract} 
\end{center}
Given $n$ samples from a population of individuals belonging to different types with unknown proportions, how do we estimate the probability of discovering a new type at the $(n+1)$-th draw? This is a classical problem in statistics, commonly referred to as the missing mass estimation problem. Recent results by Ohannessian and Dahleh \citet{Oha12} and Mossel and Ohannessian \citet{Mos15} showed: i) the impossibility of estimating (learning) the missing mass without imposing further structural assumptions on the type proportions; ii) the consistency of the Good-Turing estimator for the missing mass under the assumption that the tail of the type proportions decays to zero as a regularly varying function with parameter $\alpha\in(0,1)$. In this paper we rely on tools from Bayesian nonparametrics to provide an alternative, and simpler, proof of the impossibility of a distribution-free estimation of the missing mass. Up to our knowledge, the use of Bayesian ideas to study large sample asymptotics for the missing mass is new, and it could be of independent interest. Still relying on Bayesian nonparametric tools, we then show that under regularly varying type proportions the convergence rate of the Good-Turing estimator is the best rate that any estimator can achieve, up to a slowly varying function, and that minimax rate must be at least $n^{-\alpha/2}$. We conclude with a discussion of our results, and by conjecturing that the Good-Turing estimator is an rate optimal minimax estimator under regularly varying type proportions.
\vspace*{.2in}


\maketitle

%%%%%%%%%%%%%%%%%%%%%%%%%%%%%%%%
%%%%%%%%%%%%%%%%%%%%%%%%%%%%%%%%
%%%%%%%%%%%%%%%%%%%%%%%%%%%%%%%%
%%%%%%%%%%%%%%%%%%%%%%%%%%%%%%%%

\section{Introduction} 
Given $n$ samples from a population of individuals belonging to different types with unknown proportions, how do we estimate the probability of discovering a new type at the $(n+1)$-th draw? This is a classical problem in statistics, commonly referred to as the missing mass estimation problem. It first appeared in ecology (e.g., Fisher et al. \citet{Fis43} and Good \citet{Goo53}), and its importance has grown considerably in recent years driven by challenging applications in a wide range of scientific disciplines, such as biological and physical sciences (e.g., Kroes et al. \citet{Kro99}, Gao et al. \citet{Gao07} and Ionita-Laza et al. \citet{Ion09}), machine learning and computer science (e.g., Motwani and Vassilvitskii \citet{Mot06} and Bubeck et al. \citet{Bub13}), and information theory (e.g., Orlitsky et al. \citet{Orl04} and Ben-Hamou et al. \citet{Ben18}). To move into a concrete setting, let $P=\sum_{j\geq1}p_{j}\delta_{\theta_{j}}$ be an unknown discrete distribution, where $(\theta_{j})_{j\geq1}$ is a sequence of atoms on some measurable space and $(p_{j})_{j\geq1}$ denote the corresponding probability masses, i.e. $p_{j}\in[0,1]$ such that $\sum_{j\geq 1}p_{j}=1$. If $\mathbf{X}_{n}=(X_{1},\ldots,X_{n})$ is a collection of independent and identically distributed random variables from $P$, then we define the missing mass as
\begin{equation} \label{missmass}
M_{n}(P,\textbf{X}_{n}) = \sum_{j\geq1}p_{j} \mathbbm{1}(\theta_j \notin \textbf{X}_{n}),
 \end{equation}
where $\mathbbm{1}(\cdot)$ is the indicator function. Among various nonparametric estimators of the missing mass, both frequentist and Bayesian, the Good-Turing estimator (Good \citet{Goo53}) is arguably the most popular. It has been the subject of numerous studies, most of them in the recent years. These include, e.g., asymptotic normality and large deviations (Zhang and Zhang \citet{Zha09} and Gao \citet{Gao13}), admissibility and concentration properties (McAllester and Ortiz, \citet{McA03}, Ohannessian and Dahleh \citet{Oha12} and Ben-Hamou et al. \citet{Ben17}), consistency and convergence rates (McAllester Schapire \citet{McA00}, Wagner et al. \citet{Wag06} and Mossel and Ohannessian \citet{Mos15}), optimality and minimax properties (Orlitsky et al. \citet{Orl03} and Rajaraman et al. \citet{Raj17}).

Under the setting depicted above, let $\hat{M}_{n}(\textbf{X}_{n})$ denote an estimator of $M_{n}(P,\textbf{X}_{n})$. Motivated by the recent works of Ohannessian and Dahleh \citet{Oha12}, Mossel and Ohannessian \citet{Mos15} and Ben-Hamou et al. \citet{Ben17}, in this paper we consider the problem of consistent estimation of the missing mass under the multiplicative loss function
\begin{equation} \label{loss}
L(\hat{M}_{n}(\textbf{X}_{n}),M_{n}(P,\textbf{X}_{n}))=\left| \frac{\hat{M}_{n}(\textbf{X}_{n})}{M_{n}(P,\textbf{X}_{n})}-1 \right|.
\end{equation}
As discussed in Ohannessian and Dahleh \citet{Oha12}, the loss function \eqref{loss} is adequate for estimating small value parameters, in the sense that it allows to achieve more informative results. Such a loss function has been already used in statistics, e.g. for the estimation of small value probabilities using importance sampling (Chatterjee and Diaconis \citet{Cha18}) and for the estimation of tail probabilities in extreme value theory (Beirlant and Devroye \citet{Bei99}). Under the loss function \eqref{loss}, Ohannessian and Dahleh \citet{Oha12} showed that: i) the Good-Turing estimator may be inconsistent; ii) the Good-Turing estimator is strongly consistent if the tail of $P$ decays to zero as a regularly varying function with parameter $\alpha\in(0,1)$ (Bingham et al. \citet{Bin87}). See also Ben-Hamou et al. \citet{Ben17} for further results on missing mass estimation under regularly varying $P$. Mossel and Ohannessian \citet{Mos15} then strengthened the inconsistency result of Ohannessian and Dahleh \citet{Oha12}, showing the impossibility of estimating (learning) $M_{n}(P,\textbf{X}_{n})$ in a completely distribution-free fashion, that is without imposing further structural assumptions on $P$.

We present an alternative, and simpler, proof of the result of Mossel and Ohannessian \citet{Mos15}. Our proof relies on tools from Bayesian nonparametrics, and in particular on the use of a Dirichlet prior (Ferguson \citet{Fer73}) for the unknown distribution $P$. This allows us to exploit properties of the posterior distribution of $M_{n}(P,\textbf{X}_{n})$ to prove the impossibility of a distribution-free estimation of the missing mass, thus avoiding the winding (geometric) coupling argument of Mossel and Ohannessian \citet{Mos15}. Up to our knowledge, the use of Bayesian ideas to study large sample asymptotics for the missing mass is new, and it could be of independent interest. Motivated by the work of Ohannessian and Dahleh \citet{Oha12} and Ben-Hamou et al. \citet{Ben17} we then investigate convergence rates and minimax rates for the Good-Turing estimator under the class of $\alpha\in(0,1)$ regularly varying $P$. We still rely on tools from Bayesian nonparametrics, thus providing an original approach to tackle these problems. In particular, we make use of the two parameter Poisson-Dirichlet prior (Perman et al. \citet{Per92} and Pitman and Yor \citet{Pit97}) for the unknown distribution $P$, which is known to generate (almost surely) discrete distributions whose tail decays to zero as a regularly varying function with parameter $\alpha\in(0,1)$. See Gnedin et al. \citet{Gne07} and references therein. This allows us to exploit properties of the posterior distribution of $M_{n}(P,\textbf{X}_{n})$ to prove that: i) the convergence rate of the Good-Turing estimator is the best rate that any estimator of the missing mass can achieve, up to a slowly varying function; ii) the minimax rate must be at least $n^{-\alpha/2}$. We conclude with a discussion on the problem of deriving the minimax rate of the Good-Turing estimator, conjecturing that the Good-Turing estimator is an asymptotically optimal minimax estimator under the class of regularly varying $P$.

The paper is structured as follows. In Section \ref{Sect2} we state our main results on convergence rates and minimax rates for the Good-Turing estimator under $\alpha\in(0,1)$ regularly varying distribution $P$. Proofs of these results, as well as the alternative proof of the result of Mossel and Ohannessian \citet{Mos15}, are provided in Section \ref{Sect3}. In Section \ref{sec:discussion} we discuss open problems and possible future developments on missing mass estimation. Auxiliary results and technical lemmas are deferred to Appendix \ref{app}. The following notation is adopted throughout the paper:
\begin{itemize}
\item $[0,1]$ is the unit interval, and $\mathcal{B}([0,1])$ its Borel $\sigma$-algebra;
\item $\mathcal{P}$ is the space of discrete distributions on $[0,1]$, endowed with the smallest $\sigma$-algebra making $P\mapsto P(A)$ measurable for every $A\in \mathcal{B}([0,1])$;
\item  $P^n$ is the $n$-fold product of $P$ on $[0,1]$, and $\mathbbm{E}_{P}$ the expectation with respect $P$; for easiness of notation, we will use $\mathbbm{E}_{P}$ to denote also the expectation with respect to $P^n$.
;
\item $\ell$ is a generic slowly varying function, i.e. a  function satisfying $\ell(xc)/\ell(x)\rightarrow1$ as $x\rightarrow \infty$ for every $c>0$;
\item $C$ denotes a generic strictly positive constant that can vary in the calculations and in distinct statements;
\item Given a sequence of probabilities $(p_{j})_{j\geq1}$, $(p_{[j]})_{j\geq1}$ denotes the corresponding ordered sequence, i.e. $p_{[1]}\geq p_{[2]} \geq \ldots$;
\item Given two functions $f$ and $g$, $f\sim g$ stands for $\lim\frac{f}{g} = 1$, $f=\mathcal{O}(g)$ for $\limsup\frac{|f|}{|g|}<C$, $f=o(g)$ for $\lim\frac{f}{g}=0$;
\item $\mathcal{B}(a,b)$ is the Beta integral of parameters $a$ and $b$.
\end{itemize}

%%%%%%%%%%%%%%%%%%%%%%%%%%%%%%%%
%%%%%%%%%%%%%%%%%%%%%%%%%%%%%%%%
%%%%%%%%%%%%%%%%%%%%%%%%%%%%%%%%
%%%%%%%%%%%%%%%%%%%%%%%%%%%%%%%%

\section{Main results} \label{Sect2}
 
Let $\mathbf{X}_{n}=(X_{1},\ldots,X_{n})$ be a collection of independent and identically distributed random variables from an unknown discrete distribution $P$. The actual values taken by the observations, $X_{i}$'s, are irrelevant for the missing mass estimation problem and, without loss of generality, they can be assumed to be values in the set $[0,1]$. Therefore, $P(\cdot) = \sum_j p_j \delta_{\theta_j}(\cdot)$ is supposed to be a discrete distribution on the sample space $[0,1]$, given a sequence of atoms $\theta_{j} \in [0,1]$ and masses $p_{j}<1$ such that $\sum_{j\geq 1}p_{j}=1$. Both atoms and masses of the distribution $P$ are assumed to be unknown. Given the sample $\textbf{X}_{n}$, we are interested in estimating the missing $M_{n}(P,\textbf{X}_{n})$ defined in \eqref{missmass}, which turns out to be a  jointly measurable function of $P$ and $\textbf{X}_{n}$ as proved in Proposition \ref{prop:measurable_missing}. Given an estimator $\hat{M}_{n}(\textbf{X}_{n}):[0,1]^{n}\rightarrow[0,1]$ of  $M_{n}(P,\textbf{X}_{n})$, we will measure its statistical performance by using the multiplicative loss function defined in \eqref{loss}. As we discussed in the introduction, this loss function is suitable to study theoretical properties of parameters or functionals taking small values, and it has already been used in previous works on missing mass estimation, e.g., Ohannessian and Dahleh \citet{Oha12}, Mossel and Ohannessian \citet{Mos15} and Ben-Hamou et al. \citet{Ben17}.

A sequence of estimators $\hat{M}_{n}(\textbf{X}_{n})$ is said to be consistent for $M_{n}(P,\textbf{X}_{n})$ under parameter space $\mathcal{P}$ and loss function $L$, if the loss incurred by the estimator converges in probability to zero under all points in the parameter space. Formally, $\hat{M}_{n}(\textbf{X}_{n})$ is consistent for  $M_{n}(P,\textbf{X}_{n})$ if for all $P\in \mathcal{P}$ and for all $\epsilon>0$,
\begin{equation} \label{consistency}
P^{n}(L(\hat{M}_{n}(\textbf{X}_{n}),M_{n}(P,\textbf{X}_{n}))>\epsilon)\rightarrow 0
\end{equation}
as $n\rightarrow \infty$. Also, $\hat{M}_{n}(\textbf{X}_{n})$ is strongly consistent if \eqref{consistency} is replaced by almost sure convergence. Under this setting Mossel and Ohannessian \citet{Mos15} proved the following result.

\begin{theorem} \label{thm1}
Let $\mathcal{P}$ be the set of all discrete distributions on $[0,1]$ and $L$ be the loss function defined as \eqref{loss}. Then, there do not exist any consistent estimators for the missing mass  $M_{n}(P,\textbf{X}_{n})$, i.e. there are no estimators $\hat{M}_{n}(\textbf{X}_{n})$ satisfying \eqref{consistency}.
\end{theorem}

Mossel Ohannessian \citet{Mos15} proved Theorem \ref{thm1} by exploiting a coupling of two generalized (dithered) geometric distributions. In section Section \ref{Sect3} we present an alternative proof of Theorem \ref{thm1}. While the proof of Mossel and Ohannessian \citet{Mos15} has the merit to be constructive, our approach has the merit to be simpler and it provides a new way to face these type of problems, which mainly relies on Bayesian nonparametric techniques. Similar Bayesian nonparametric arguments will then be crucial in order to study of convergence rates and minimax rates of the Good-Turing estimator under the class of $\alpha\in(0,1)$ regularly varying $P$.

Roughly speaking, Theorem \ref{thm1} proves that  any asymptotic result holding uniformly over a set of possible distributions, the parametric space $\mathcal{P}$ must be restricted to a suitable subclass. In particular, from the proof of Theorem \ref{thm1} we see that some conditions have to be imposed on the tail decay of the elements of the parameter space. That is, from the proof of Theorem \ref{thm1} we deduce that there are no consistent estimators for the class of distributions sampled from a Dirichlet process. From Kingman \citet{Kin75} (Equation 65), we have that, if $P$ is sampled from a Dirichlet process, its sequence of ordered masses behaves like $\log p_{[j]}  \sim -jC$, as $j\rightarrow +\infty$. Therefore, the tail of $P$ has approximately exponential form, resembling a geometric distribution and satisfying $p_{[j]}=o(j^{-\frac{1}{\alpha}})$ for every $\alpha\in (0,1)$. Indeed, a geometric distribution was used in Ohannessian and Dahleh \citet{Oha12} as an example to prove that the Good-Turing estimator can be inconsistent. Theorem \ref{thm1} shows that, under this very light regime, any estimator of the missing mass, not just the Good-Turing, fails to be consistent under multiplicative loss. This motivates us to consider the class of $P$s having heavy enough tails. This will be the subject of the rest of this section.

\subsection{Consistency under regularly varying $P$} \label{sec:rv}

In this section we recall the Good-Turing estimator (Good \citet{Goo53}) of the missing mass $M_{n}(P,\textbf{X}_{n})$, and we study its convergence rate and minimax risk for regularly varying $P$. The definition of the Good-Turing estimator makes use of  the proportion of unique values in the sample to estimate the missing mass. Let $Y_{n,j}(\textbf{X}_{n})$ be the number of times the value $\theta_{j}$ is observed in the sample $\textbf{X}_{n}$, i.e.,
\begin{equation*}
Y_{n,j}(\textbf{X}_{n}) = \sum_{i=1}^{n} \mathbbm{1}(X_i = \theta_j).
\end{equation*}
Furthermore, let $K_{n,r}(\textbf{X}_{n})$ and $K_{n}(\textbf{X}_{n})$ be the number of values observed $1\leq r\leq n$ times and the total number of distinct values, respectively, observed in $\textbf{X}_{n}$, i.e., 
\begin{equation*}
K_{n,r}(\textbf{X}_{n}) = \sum_{j=1}^{\infty} \mathbbm{1}(Y_{n,j} = r)\qquad K_{n}(\textbf{X}_{n}) =\sum_{r=1}^{n}K_{n,r}(\textbf{X}_{n}). 
\end{equation*}
The Good-Turing estimator of $M_{n}(P,\textbf{X}_{n})$ is defined in terms of the statistic $K_{n,1}(\textbf{X}_{n})$, that is 
\begin{equation} \label{GoodTuring}
\hat{GT}(\textbf{X}_{n})=\frac{K_{n,1}(\textbf{X}_{n})}{n}.
\end{equation} 
Ohannessian and Dahleh \citet{Oha12} first showed the inconsistency of $\hat{GT}(\textbf{X}_{n})$ under the choice of $P$ being a geometric distribution. In the same paper, it is shown that under the assumption that the tail of $P$ decays to zero as a regularly varying function with parameter $\alpha\in(0,1)$, the Good-Turing estimator is strongly consistent. This latter result was generalized to the range $\alpha\in (0,1]$ in Ben-Hamou et al. \citet{Ben17}.

The assumption of regularly varying $P$ is a generalization of the power law tail decay, adding some more flexibility by the introduction of the slowly varying function $\ell$. Power-law distributions are observed in the empirical distributions  of many quantities in different applied areas, and their study have attracted a lot of interest in recent years. For extensive discussions of power laws in empirical data and their properties, the reader is referred to Mitzenmacher \citet{Mit04}, Goldwater et al. \citet{Gol06}, Newman \citet{New03}, Clauset et al. \citet{Cla09} and Sornette \citet{Sor06}. Restricting the parameter space to probability distributions having regularly varying tail is not a mere technical assumption and, on the contrary, it represents a natural subset of the parameter space to consider, which we expect to contain the true data generating distribution for many different applications.

To move into the concrete setting of regular variation (Bingham et al. \citet{Bin87}), for every $P\in \mathcal{P}$ we define a counting measure on $[0,1]$ as $\nu_{P} (dx)=\sum_j \delta_{p_j}(dx)$, with corresponding tail function defined as $\vec{\nu}_{P}(x) = \nu ([x,+\infty))$ for all $x>0$. Then a distribution $P\in \mathcal{P}$ is said to be regularly varying with parameter $\alpha\in (0,1)$ if
\begin{equation} \label{regular variation}
\vec{\nu}_{P}(x) \stackrel{x\downarrow0}{\sim} x^{-\alpha} \ell(1/x),
\end{equation}
where $\ell$ is a slowly varying function. From Lemma 22 and Proposition 23 of Gnedin et al. \citet{Gne07}, \eqref{regular variation} is equivalent to the more explicit condition in term of ordered masses of $P$
\begin{equation}\label{reg_var_p}
p_{[j]} \stackrel{j\uparrow \infty}{\sim} j^{-1/\alpha} \ell_*(j),
\end{equation}
where $\ell_*$ is a slowly varying function depending on $\ell$. We denote by $\mathcal{P}_{RV_\alpha}\subseteq\mathcal{P}$ the set of all regularly varying distribution on $[0,1]$ with parameter $\alpha$. From \eqref{reg_var_p} it is clear that such a class includes distributions having power-law tail decay, which correspond to the particular case of $\ell_*$ being a constant, which is equivalent to $\ell$ being constant. We denote the class of distributions having power law tail decay by $\mathcal{P}_{PL_\alpha}\subseteq\mathcal{P}_{RV_\alpha}$. In the following results, we will restrict our attention to the estimation problem under restricted parameter spaces $\mathcal{P}_{RV_\alpha}$ and $\mathcal{P}_{PL_\alpha}$. 

From Ohannessian and Dahleh \citet{Oha12} it is known that the Good-Turing estimator is consistent under all $P$s belonging to the regularly varying class $\mathcal{P}_{RV_\alpha}$. Focusing attention on the class $\mathcal{P}_{RV_\alpha}$, in the next proposition we refine the result of Ohannessian and Dahleh \citet{Oha12} by studying the rate at which the multiplicative loss of the Good-Turing estimator converges to zero. A sequence $(r_{n})_{n\in\mathbbm{N}}$ is a convergence rate of an estimator $\hat{M}_{n}(\textbf{X}_{n})$ for the distribution $P\in \mathcal{P}$ if
\begin{equation*}
\lim_n P^{n}(L(\hat{M}_{n}(\textbf{X}_{n}),M_{n}(P,\textbf{X}_{n})) > T_{n} r_n) = 0
\end{equation*}
for all sequences $T_{n}\rightarrow \infty$. The next proposition shows the rate of convergence of the Good-Turing estimator $\hat{GT}(\textbf{X}_{n})$ is $n^{-\alpha /2}\ell^{-1/2}(n)$. The proof is omitted because it follows from Proposition \ref{prop2} below along with a simple application of Markov's inequality. 

\begin{proposition} \label{prop1}
Let $\hat{GT}(\textbf{X}_{n})$ be the Good-Turing estimator, defined in \eqref{GoodTuring}. Then, for every $P\in \mathcal{P}_{RV_\alpha}$ and for all $T_{n}\rightarrow\infty$,
\begin{equation} \label{rateGT}
\lim_n P^{n}(L(\hat{GT}(\textbf{X}_{n}),M_{n}(P,\textbf{X}_{n})) > T_{n} n^{-\alpha/2}\ell^{-1/2}(n)) = 0,
\end{equation}
where $\ell$ in \eqref{rateGT} is the slowly varying function specific to $P$ appearing in \eqref{regular variation}.
Therefore, up to slowly varying functions, $n^{-\alpha/2}$ is a convergence rate for the Good-Turing estimator within the class $\mathcal{P}_{RV_\alpha}$.
\end{proposition}

As a further result on convergence rate of the Good-Turing estimator, in Theorem \ref{thm3} we show that the convergence rate achieved by the Good-Turing estimator is actually almost the best convergence rate any estimator of $M_{n}(\textbf{X}_{n},P)$ can achieve. Specifically, for any other estimator, it is possible to find a point $P\in \mathcal{P}_{PL_\alpha}$ for which the rate of convergence is not faster than $n^{-\alpha /2}$. 

\begin{theorem} \label{thm3}
For any estimator $\hat{M}_{n}(\textbf{X}_{n})$, there exists $P\in \mathcal{P}_{PL_\alpha} \subset \mathcal{P}_{RV_\alpha}$  such that for every $T_{n}\rightarrow 0$
\begin{equation}\label{bananas1}
 \liminf_n P^{n} (L(\hat{M}_{n}(\textbf{X}_{n}),M_{n}(P,\textbf{X}_{n})) < T_{n} n^{-\alpha/2}) = 0.
\end{equation}
Therefore the convergence rate of $\hat{M}_{n}(\textbf{X}_{n})$ cannot be faster than $n^{-\alpha/2}$.
\end{theorem}

Proposition \ref{prop1} and Theorem \ref{thm3} together show that the Good-Turing estimator achieves the best convergence rate up to possibly a slowly varying function. In particular, if the distribution $P$ has a power-law decay, i.e. $P\in \mathcal{P}_{PL_\alpha}$, the two rates match and the Good-Turing estimator achieves the best rate possible. In particular, because $\hat{GT}(\mathbf{X}_{n})$ does not depend on $\alpha$, it follows that the Good-Turing estimator is actually \textit{rate adaptive} for the class of power law distributions, $\mathcal{P}_{PL}=\cup_{0<\alpha<1}\mathcal{P}_{PL_{\alpha}}$ However, for a general $P\in \mathcal{P}_{RV_\alpha}$ we do not know whether the two rates of Proposition \ref{prop1} and Theorem \ref{thm3} may be improved to make them match or they are not.

As a final result, in the next theorem we consider the asymptotic minimax estimation risk for the missing mass under the loss function \eqref{loss} and with parameter space $\mathcal{P}_{PL_\alpha}$. Theorem \ref{thm4} provides with a lower bound for the estimation risk of this statistical problem, showing that the minimax rate is not smaller than $n^{-\alpha/2}$. 

\begin{theorem} \label{thm4}
Let $\mathcal{P}_{PL_\alpha}$ be the class of discrete distributions on $[0,1]$ with power law tail function and let $L$ denote the multiplicative loss function \eqref{loss}. Then, there exists a positive constant $C>0$ such that
\begin{equation*}
\liminf_{n} n^{\alpha/2} \inf_{\hat{M}_{n}(\textbf{X}_{n})}\sup_{P\in \mathcal{P}_{PL_\alpha}}\mathbbm{E}_{P}\left( L(\hat{M}_{n}(\textbf{X}_{n}),M_{n}(P,\textbf{X}_{n})) \right)>C
\end{equation*}
where the infimum is taken over all possible estimators $\hat{M}_{n}(\textbf{X}_{n})$.
\end{theorem} 

The lower bound of Theorem \ref{thm4} can be used to derive the minimax rate, by matching it with appropriate upper bounds of specific estimators of the missing mass. This lower bound trivially still holds for any parametric set larger than  $\mathcal{P}_{PL_\alpha}$ and, therefore, the theorem also provides with a lower bound of the estimation risk under the larger parameter space $\mathcal{P}_{RV_\alpha}$ . In the next Proposition, we show that for a fixed distribution $P\in\mathcal{P}_{RV_\alpha}$, the Good-Turing estimator achieves the best possible rate of Theorem \ref{thm4} up to a slowly varying term. 

\begin{proposition} \label{prop2}
Let $\hat{GT}(\textbf{X}_{n})$ be the Good-Turing estimator and let $P\in \mathcal{P}_{RV_\alpha}$. Then, there exists a finite constant $C$ such that for every $n$
\begin{equation} \label{riskGT}
\mathbbm{E}_{P}(L(\hat{GT}(\textbf{X}_{n}),M_{n}(P,\textbf{X}_{n}))) \leq  C n^{-\alpha/2}\ell^{-1/2}(n),
\end{equation}
where $\ell$ is the slowly varying function specific to $P$ appearing in \eqref{regular variation}.
\end{proposition}

Extending Proposition \ref{prop2} to hold uniformly over $\mathcal{P}_{RV_\alpha}$ is an open problem and probably requires a careful control over the size of $\mathcal{P}_{RV_\alpha}$. Indeed, the classes of distributions we are considering are defined through the asymptotic properties of their elements, while to obtain minimax results we need a control for each $n \in \mathbbm{N}$. Even though Proposition \ref{prop2} does not directly provide with the minimax rate of the Good-Turing estimator, it still provides with a sanity check for its asymptotic risk. Specifically, Proposition \ref{prop2} implies that for every $P \in \mathcal{P}_{PL_\alpha}$, 
\begin{displaymath}
\limsup_n n^{\alpha/2}\mathbbm{E}_{P}(L(\hat{GT}(\textbf{X}_{n}),M_{n}(P,\textbf{X}_{n}))) < +\infty
\end{displaymath}
Moreover, from a minor change at the beginning of the proof of Theorem \ref{thm4}, we can also prove that for every estimator $\hat{M}_{n}(\textbf{X}_{n})$ and every sequence $(T_n)_n$ diverging to infinity, we can find an element $P \in \mathcal{P}_{PL_\alpha} $ such that  $\limsup_n T_n n^{\alpha/2}\mathbbm{E}_{P}(L(\hat{M}_{n}(\textbf{X}_{n}),M_{n}(P,\textbf{X}_{n}))) = +\infty $. This leads us to conjecture that the Good-Turing estimator should be a rate optimal minimax estimator.

%%%%%%%%%%%%%%%%%%%%%%%%%%%%%%%%
%%%%%%%%%%%%%%%%%%%%%%%%%%%%%%%%
%%%%%%%%%%%%%%%%%%%%%%%%%%%%%%%%
%%%%%%%%%%%%%%%%%%%%%%%%%%%%%%%%

\section{Discussion} \label{sec:discussion}
In this paper we have considered the problem of consistent estimation of the missing mass under a suitable multiplicative loss function. We have presented an alternative, and simpler, proof of the result by Mossel and Ohannessian \citet{Mos15} on the impossibility of a distribution-free estimation of the missing mass. Our results relies on novel arguments from Bayesian nonparametric statistics, which are then exploited to study convergence rates and minimax rates of the Good-Turing estimator under the class of $\alpha\in(0,1)$ regularly varying $P$. In Proposition \ref{prop1} and Theorem \ref{thm3} it has been shown that, within the class $\mathcal{P}_{PL_\alpha}$, the Good-Turing estimator achieves the best convergence rate possible, while for the  class, $\mathcal{P}_{RV_\alpha}$, this rate is the best up to a slowly varying function. An open problem is to understand weather this additional slowly varying term is intrinsic to the problem or our results can actually be improved to make the rate of the Good-Turing estimator matches the best possible rate also within the class of regularly varying distributions. Under the restricted parametric spaces, in Theorem \ref{thm4} we have provided a lower bound for the asymptotic risk. This bound can be used to compare estimators from a minimax point of view, by finding suitable upper bounds matching the lower bound rate. In particular, in Proposition \ref{prop2} we have shown that the asymptotic rate of the risk of the Good-Turing estimator matches the lower bound rate, up to a slowly varying function. However, the rate of Proposition \ref{prop2} is a pointwise result, for a fixed $P\in \mathcal{P}_{RV_{\alpha}}$. An open problem is to extend Proposition \ref{prop2} to the uniform case, when considering  the supremum of the risk over all $P\in \mathcal{P}_{RV_{\alpha}}$. This extension probably requires a careful analysis and control of the size of this parameter space $\mathcal{P}_{RV_{\alpha}}$. Work on this is ongoing. 

%%%%%%%%%%%%%%%%%%%%%%%%%%%%%%%%
%%%%%%%%%%%%%%%%%%%%%%%%%%%%%%%%
%%%%%%%%%%%%%%%%%%%%%%%%%%%%%%%%
%%%%%%%%%%%%%%%%%%%%%%%%%%%%%%%%

\section{Proofs} \label{Sect3}
In this section we will prove all the theorems stated in Section \ref{Sect2}. The proofs of some technical auxiliary results are postponed to Appendix \ref{app}. We start with a simple lemma that will be useful in the sequel.
\begin{lemma}\label{lemma_inv_L}
For $\epsilon < 1/2$ and $a,b \geq 0$, $L(a,b) \leq \epsilon$ implies $L(b,a) \leq 2\epsilon.$
\end{lemma}
\textsc{proof}
Let $a,b$ be positive real numbers, $\epsilon < 1/2$ and suppose $L(a,b)=|\frac{a}{b} - 1| \leq \epsilon$. Straightforwardly, have that
\begin{equation}
-b\epsilon  \leq a-b \leq b\epsilon \label{ineq:inv_d1}
\end{equation}
From the lower bound of \eqref{ineq:inv_d1}, $a \geq (1-\epsilon) b \geq b/2$ and, therefore,
$ \frac{1}{a} \leq \frac{2}{b}$. Multiplying (\ref{ineq:inv_d1}) by this last inequality, we have 
$L(b,a) \leq 2\epsilon $.
\qed

\subsection{Proof of Theorem \ref{thm1}}
We are going to show that for every estimator $\hat{M}_{n}(\textbf{X}_{n})$, there exists $\epsilon>0$ such that
\begin{equation}
 \sup_{P \in \mathcal{P}} \limsup_n P^{n}(L(\hat{M}_{n}(\textbf{X}_{n}),M_{n}(P,\textbf{X}_{n}))>\epsilon) > 0 \label{bananaquatre}
\end{equation}
and, therefore, there exists $P \in \mathcal{P}$ such that $L(\hat{M}_{n}(\textbf{X}_{n}),M_{n}(P,\textbf{X}_{n}))$ does not converge to zero in probability. 

First, note that for $\epsilon<1/2$, Lemma (\ref{lemma_inv_L}) implies that
\begin{equation} \label{bananabis}
P^{n}(L(\hat{M}_{n}(\textbf{X}_{n}),M_{n}(P,\textbf{X}_{n}))>\epsilon)  \geq P^{n}(L(M_{n}(P,\textbf{X}_{n}),\hat{M}_{n}(\textbf{X}_{n}))> 2\epsilon ), 
\end{equation}
so it is sufficient to show that there exists $0<\epsilon<1$ such that for every estimator $\hat{M}_{n}(\textbf{X}_{n})$,
\begin{equation} \label{bananatris}
 \sup_{P \in \mathcal{P}} \limsup_n P^{n}(L(M_{n}(P,\textbf{X}_{n}),\hat{M}_{n}(\textbf{X}_{n}))>\epsilon) > 0.
\end{equation}
We will prove that \eqref{bananatris} holds for all $0<\epsilon<1/4$ (and therefore (\ref{bananaquatre}) holds for any $0<\epsilon<1/8$). Let $\epsilon \in (0,1/4)$. Let $\text{DP}_{\gamma}$ denote the Dirichlet process measure on $\mathcal{P}$ (Ferguson \citet{Fer73}), with base measure $\gamma$ on $[0,1]$. We choose $\gamma$ uniform, i.e. $\gamma(d\theta)=\mathbbm{1}(0<d\theta<1)$. Now, we can lower bound the supremum in \eqref{bananatris} by an average over $\mathcal{P}$ with respect to $\text{DP}_{\gamma}$ and then swap the integration by Fubini theorem, therefore 

\begin{align*}
& \sup_{P \in \mathcal{P}}  \limsup_n P^{n}(L(M_{n}(P,\textbf{X}_{n}),\hat{M}_{n}(\textbf{X}_{n}))>\epsilon)  \\
 & \  \geq \int_{\mathcal{P}} \limsup_n P^{n}(L(M_{n}(P,\textbf{X}_{n}),\hat{M}_{n}(\textbf{X}_{n}))>\epsilon) \text{DP}_{\gamma}(dP)  \\
 & \ \geq \limsup_n \int_{\mathcal{P}} \int_{[0,1]^{n}} \mathbbm{1}(L(M_{n}(P,\textbf{X}_{n}),\hat{M}_{n}(\textbf{X}_{n}))>\epsilon) P^{n}(d\textbf{X}_{n}) \text{DP}_{\gamma}(dP)  \\
  & \ = \limsup_n \int_{[0,1]^{n}} \int_{\mathcal{P}} \mathbbm{1}(L(M_{n}(P,\textbf{X}_{n}),\hat{M}_{n}(\textbf{X}_{n}))>\epsilon) \text{DP}_{\gamma + \sum_{i=1}^{n}\delta_{X_{i}}}(dP) P^{n}_{\text{DP}_{\gamma}}(d\textbf{X}_{n})    \\
    & \ \geq \limsup_n \int_{[0,1]^{n}} \inf_{x\geq 0}  \int_{\mathcal{P}} \mathbbm{1}(L(M_{n}(P,\textbf{X}_{n}),x)>\epsilon) \text{DP}_{\gamma + \sum_{i=1}^{n}\delta_{X_{i}}}(dP) d\textbf{X}_{n}   
\end{align*}
where the first inequality follows since we can lower bound the supremum by an average, the second from reverse Fatou's lemma, the equality comes by swapping the marginal of $P$ and conditional of $\textbf{X}_{n}$ given $P$ with the marginal of $\textbf{X}_{n}$, denoted $P^{n}_{\text{DP}_{\gamma}}$, and the conditional of $P$ given $\textbf{X}_{n}$, the last inequality follows since we are considering the infimum over all possible values of $\hat{M}_{n}(\textbf{X}_{n})$. Also recall that, when $P$ is distributed as $\text{DP}_{\gamma}$, then the marginal of $\textbf{X}_{n}$, $P^{n}_{\text{DP}_{\gamma}}$, is a Generalized Polya urn, while the conditional of $P$ given $\textbf{X}_{n}$ is  $\text{DP}_{\gamma + \sum_{i=1}^{n}\delta_{X_{i}}}$ (see Theorem 4.6 and subsection 4.1.4 of Ghosal and Van der Vaart \citet{Gho17}).
 
From Proposition \ref{propPitman} (Appendix \ref{app}),  $M_{n}(P,\textbf{X}_{n})$ under the posterior distribution $\text{DP}_{\gamma + \sum_{i=1}^{n}\delta_{X_{i}}}$ is distributed according to a Beta random variable $\text{Beta}(1,n)$. Therefore,
\begin{equation} \label{banana8}
 \int_{\mathcal{P}} \mathbbm{1}(L(M_{n}(P,\textbf{X}_{n}),x)>\epsilon) \text{DP}_{\gamma + \sum_{i=1}^{n}\delta_{X_{i}}}(dP)= \mathbbm{P}\left(\left| \frac{Z}{x}-1 \right| >\epsilon \right)
\end{equation}
where $Z\sim \text{Beta}(1,n)$. We are now going to lower bound the probability of the event on the right hand side of \eqref{banana8}. 
First let us consider $x \in (0,\frac{1}{1+\epsilon}]$ and $n\geq 2$
\begin{align}
\mathbbm{P}\left(\left| \frac{Z}{x}-1 \right| >\epsilon \right) &= \mathbbm{P} (Z > (1+\epsilon)x) + \mathbbm{P} (Z < (1-\epsilon)x) \nonumber \\
& = 1 + (1-(1+\epsilon)x)^n - (1-(1-\epsilon)x)^n \nonumber \\
&= 1 - 2 x \epsilon \sum\limits_{k=0}^{n-1} (1-(1+\epsilon)x)^{n-1-k} (1-(1-\epsilon)x)^k \label{formula1} \\
&\geq 1 - 2 x \epsilon \sum\limits_{k=0}^{n-1} (1-(1-\epsilon)x)^{n-1} \nonumber \\
&= 1 - 2 x \epsilon n(1-(1-\epsilon)x)^{n-1} \nonumber  \\
&\geq  1 - 2 \frac{\epsilon}{(1-\epsilon)} n (1-\epsilon) x (1-(1-\epsilon)x)^{n-1} \nonumber \\
&\geq  1 - 2\frac{\epsilon}{(1-\epsilon)} (1-1/n)^{n-1} \label{formula2} \\
&\geq  1 - \frac{2\epsilon}{(1-\epsilon)} \nonumber
\end{align}
where we have used $a^n-b^n = (a-b) \sum\limits_{k=0}^{n-1} a^{n-1-k} b^k$ in \eqref{formula1} and that the maximum of the function $x \mapsto x(1-x)^{n-1}$ is achieved in $1/n$ in \eqref{formula2}. Now let $x > \frac{1}{1+\epsilon}$, noticing that $\frac{2\epsilon}{(1+\epsilon)} < 1$, it comes that
\begin{align*}
\mathbbm{P}\left(\left| \frac{Z}{x}-1 \right| >\epsilon \right) &= \mathbbm{P} \left( Z < (1-\epsilon)x\right)\geq \mathbbm{P} \left( Z < \frac{1-\epsilon}{1+\epsilon}\right) = 1 - 2^n \frac{\epsilon^n}{(1+\epsilon)^n}\\
&\geq 1 - \frac{2\epsilon}{(1+\epsilon)}
\geq 1 - \frac{2\epsilon}{(1-\epsilon)}.
\end{align*}
Therefore, $\mathbbm{P}\left(\left| \frac{Z}{x}-1 \right| >\epsilon \right) \geq 1 - \frac{2\epsilon}{(1-\epsilon)}$ for all $x\geq 0$ and $n\geq2$. Plugging this estimate in place of \eqref{banana8}, we obtain
\begin{equation*}
 \sup_{P \in \mathcal{P}}  \limsup_n P^{n}(L(M_{n}(P,\textbf{X}_{n}),\hat{M}_{n}(\textbf{X}_{n}))>\epsilon) \geq  1 - \frac{2\epsilon}{(1-\epsilon)}
\end{equation*}
and the right hand side is strictly positive for all $0<\epsilon<1/4$.

\subsection{Proof of Theorem \ref{thm3}}
Let $(T_n)_n$ any non-negative sequence converging to $0$. We will show that for any estimator $\hat{M}_{n}(\textbf{X}_{n})$, 
\begin{equation} \label{ciao}
\inf_{P\in \mathcal{P}_{RV_\alpha}} \liminf_n P^{n}(L(\hat{M}_{n}(\textbf{X}_{n}),M_{n}(P,\textbf{X}_{n})) < T_n n^{-\alpha/2}) = 0.
\end{equation}

Let us denote by $SP_{\alpha}$ the law of a stable process on $[0,1]$ of parameter $\alpha$. This a subordinator with Levy intensity, $\nu (d\omega) = \frac{\alpha}{\Gamma(1-\alpha)}\omega ^{-1-\alpha} d\omega$. See Kingman \citet{Kin75}, Lijoi and Pr\"unster \citet{Lij10} and Pitman \citet{Pit06} for details and additional references. Because 
of $\nu ([x,\infty))= \frac{x^{-\alpha}}{\Gamma(1-\alpha)}$, the stable process samples probability measures belonging to $\mathcal{P}_{PL_\alpha}$. Now we can upper bound the infimum in \eqref{ciao} by an average with respect to $SP_{\alpha}$,
\begin{align*}
& \inf_{P\in \mathcal{P}_{RV_\alpha}} \ \liminf_n P^{n}(L(\hat{M}_{n}(\textbf{X}_{n}),M_{n}(P,\textbf{X}_{n})) < T_n n^{-\alpha/2})  \\
& \ \ \leq \int_{\mathcal{P}}  \liminf_n P^{n}(L(\hat{M}_{n}(\textbf{X}_{n}),M_{n}(P,\textbf{X}_{n})) < T_n n^{-\alpha/2}) SP_{\alpha}(dP) \\
& \ \ \leq  \liminf_n \int_{\mathcal{P}} \int_{[0,1]^{n}} \mathbbm{1}(L(\hat{M}_{n}(\textbf{X}_{n}),M_{n}(P,\textbf{X}_{n})) < T_n n^{-\alpha/2}) P^{n}(d\textbf{X}_{n})  SP_{\alpha}(dP)
\end{align*}
where the last equality follows by applying Fatou's Lemma.

Take $n$ large enough so that $T_n n^{-\alpha/2} < 1/2$. Let us denote by $P^{n}_{SP_{\alpha}}$ the marginal law of the observations under an $\alpha$-stable process, when $P$ is integrated out, i.e. the probability measure on $[0,1]^{n}$ defined as
$P^{n}_{SP_{\alpha}}(A)=\int_{\mathcal{P}}P^{n}(A)SP_{\alpha}(dP)$ for all $A\in\mathcal{B}([0,1]^{n})$
. We swap the integration of the marginal of $P$ and the conditional of $\textbf{X}_{n}$ given $P$ with the marginal of $\textbf{X}_{n}$ and the conditional of $P$ given $\textbf{X}_{n}$ and then apply Lemma \ref{lemma_inv_L} to obtain 

\begin{align*}
& \int_{\mathcal{P}} \int_{[0,1]^{n}}  \mathbbm{1}(L(\hat{M}_{n}(\textbf{X}_{n}),M_{n}(P,\textbf{X}_{n})) < T_n n^{-\alpha/2}) P^{n}(d\textbf{X}_{n})  SP_{\alpha}(dP) \\
& = \int_{[0,1]^{n}} \int_{\mathcal{P}}  \mathbbm{1}(L(\hat{M}_{n}(\textbf{X}_{n}),M_{n}(P,\textbf{X}_{n})) < T_n n^{-\alpha/2})  SP_{\alpha}|\textbf{X}_{n}(dP)  P^{n}_{SP_{\alpha}}(d\textbf{X}_{n}) \\
& \leq \int_{[0,1]^{n}} \int_{\mathcal{P}}  \mathbbm{1}(L(M_{n}(P,\textbf{X}_{n}),\hat{M}_{n}(\textbf{X}_{n}))) < 2 T_n n^{-\alpha/2})  SP_{\alpha}|\textbf{X}_{n}(dP)  P^{n}_{SP_{\alpha}}(d\textbf{X}_{n})
\end{align*}

where $SP_{\alpha}|\textbf{X}_{n}$ denotes the posterior of $P$ given the sample.
Therefore, taking $s > 1$, we can upper bound the quantity appearing on the l.h.s. of \eqref{ciao} by

\begin{equation} \label{ciaociao}
\begin{split}
& \limsup_n  \int_{[0,1]^{n}} \mathbbm{1}(K_n(\textbf{X}_{n}) \in (n^\alpha/s , s n^\alpha))     \\
& \times \int_{\mathcal{P}} \mathbbm{1}(L(M_{n}(P,\textbf{X}_{n}),\hat{M}_{n}(\textbf{X}_{n}))) < 2T_n n^{-\alpha/2})   SP_{\alpha}|\textbf{X}_{n}(dP)  P^{n}_{SP_{\alpha}}(d\textbf{X}_{n}) \\
 & \qquad\qquad\qquad+   \limsup_n  \  P^{n}_{SP_{\alpha}} ( K_n(\textbf{X}_{n}) \not \in (n^\alpha/s, s n^\alpha) ) 
 \end{split}
\end{equation}  

We will upper-bound the two terms of the sum in \eqref{ciaociao} independently.
Let us focus on the first term of \eqref{ciaociao}. Let $n$ large enough so that $  3 < \frac{\alpha n^\alpha}{s} < \alpha s n^\alpha < n - 3 $ and $T_n n^{-\alpha/2} < 1/4$. From Proposition \ref{propPitman}, under the posterior $SP_{\alpha}|\textbf{X}_{n}$, $M_{n}(P,\textbf{X}_{n})$ is distributed according to a Beta random variable $\text{Beta}(\alpha K_{n}(\textbf{X}_{n}),n-\alpha K_{n}(\textbf{X}_{n}))$. Let us denote $a(\textbf{X}_{n})=\alpha K_{n}(\textbf{X}_{n})$, $b(\textbf{X}_{n})=n-\alpha K_{n}(\textbf{X}_{n})$, and for easiness of notation we will simply write $a$ and $b$ in the following calculations. Also let $F_{a,b}$ be the cumulative distribution function of the beta random variable $\text{Beta}(a,b)$. From Proposition \ref{propPitman} in the Appendix, we have that 
\begin{align*}
&\int_{\mathcal{P}} \mathbbm{1}(L(M_{n}(P,\textbf{X}_{n}),\hat{M}_{n}(\textbf{X}_{n}))) < 2T_n n^{-\alpha/2})  SP_{\alpha}|\textbf{X}_{n}(dP)   \\
& \qquad =  F_{a,b}( (1+2 T_n n^{-\alpha/2}) \hat{M}_{n}(\textbf{X}_{n}) ) - F_{a,b}( (1-2 T_n n^{-\alpha/2}) \hat{M}_{n}(\textbf{X}_{n}) ) \\
 & \qquad\leq\sup_{x\in [0,1]} \left( F_{a,b}( (1+2 T_n n^{-\alpha/2}) x) - F_{a,b}( (1-2 T_n n^{-\alpha/2}) x) \right).
\end{align*}
Consider the function $\psi : \mathbbm{R}_+ \rightarrow [0,1]$ defined by
$$ \psi(x) = F_{a,b}( (1+2 T_n n^{-\alpha/2}) x) - F_{a,b}( (1-2 T_n n^{-\alpha/2}) x). $$
Notice that $\psi \in \mathcal{C}^2$ and that $\psi(0) = \lim_{x\rightarrow\infty} \psi(x) = 0$. Therefore, $\psi$ reaches its maximum in $x^*(a,b) \in \mathbbm{R}_+$ (denoted $x^*$ for easiness of notation) satisfying
\begin{align*}
 \psi'(x^*) = (1+ & 2 T_n n^{-\alpha/2}) f_{a,b}( (1+2 T_n n^{-\alpha/2}) x^*)  \\
& - (1-2 T_n n^{-\alpha/2}) f_{a,b}( (1-2 T_n n^{-\alpha/2}) x^*)  = 0,
\end{align*}
where $f_{a,b}$ denotes the density function of the $\text{Beta}(a,b)$ distribution. On the event $K_n(\textbf{X}_{n}) \in (n^\alpha/s, s n^\alpha)$, we have $a,b > 3$, and so $f_{a,b}$ is bell-shaped with second inflexion point, 
\begin{align*}
\kappa(a,b) &=  \frac{a-1}{a+b-2} + \frac{ \sqrt{ \frac{(a-1)(b-1)}{a+b-3} } }{a+b-2} \\
& \leq \frac{\alpha s n^\alpha}{n-2} +  \frac{\sqrt{\alpha s n^\alpha}}{n-2} \\
& \leq \frac{2 \alpha s n^\alpha}{n-2}.
\end{align*}
Therefore, $f'_{a,b}$ is non decreasing on the interval $[\kappa(a,b),\infty)$ and as a consequence, $ \psi''$ is non negative on $[\frac{\kappa(a,b)}{(1-2 T_n n^{-\alpha/2})},\infty)$, from which we can deduce that $\psi'$ is non decreasing on the same interval. Now, since $\lim_{x\rightarrow\infty} \psi'(x) = 0$, it follows that $\psi'(x) \leq 0$ on $[\frac{\kappa(a,b)}{(1-2 T_n n^{-\alpha/2})},\infty)$. Therefore,
$$ x^*(a,b) \leq \frac{\kappa(a,b)}{(1-2 T_n n^{-\alpha/2})} \leq \frac{2 \alpha s n^\alpha}{(n-2)(1-2 T_n n^{-\alpha/2})} \leq \frac{4 \alpha s n^\alpha}{n-2}. $$
We can now upper-bound $\sup_{x\geq 0} \psi(x)$ as follows:
\begin{align*}
\sup_{x\geq 0} \psi(x) &= \psi(x^*) \leq 4 T_n n^{-\alpha/2} x^* \sup_{x\geq 0} f_{a,b}(x) 
\leq 16 T_n n^{-\alpha/2} \frac{\alpha s n^\alpha}{n-2} \sup_{x\geq 0} f_{a,b}(x)
\end{align*}
From Lemma \ref{lemma1} in Appendix \ref{app}, it follows that, on the event $K_n \in (n^\alpha/s,  s n^\alpha)$, for $n$ large enough,
\begin{align*}
\sup_{x\geq 0} \psi(x) &\leq 128 T_n n^{-\alpha/2} \frac{\alpha s n^\alpha}{n-2} (a+b)^{3/2} a^{-1/2} b^{-1/2} \\
&\leq  128 T_n n^{-\alpha/2} \frac{\alpha s n^\alpha}{n-2} n^{3/2} (\alpha n^\alpha/s )^{-1/2} (n-s \alpha n^\alpha)^{-1/2} = T_n g(\alpha,s,n)
\end{align*}
and notice that $\limsup\limits_{n\rightarrow +\infty} T_n g(\alpha,s,n) = 0$.\\

From all previous computations, we can deduce that, on the event $K_n(\textbf{X}_{n}) \in (1/s n^\alpha, s n^\alpha)$, we can find $n_0 (\alpha,s)$ ($n_0$ does not depend on the value of $K_n$) such that for all $n \geq n_0(\alpha,s)$, the inequality
\begin{equation*}
\int_{\mathcal{P}}  \mathbbm{1}(L(M_{n}(P,\textbf{X}_{n}),\hat{M}_{n}(\textbf{X}_{n}))) < 2T_n n^{-\alpha/2})  SP_{\alpha}|\textbf{X}_{n}(dP)  \leq T_n g(\alpha,s,n) 
\end{equation*}
holds, leading to
\begin{align*}
& \limsup_n  \int_{[0,1]^{n}} \mathbbm{1}(K_n(\textbf{X}_{n}) \in (1/s n^\alpha, s n^\alpha))  \nonumber \\
 & \qquad \times \int_{\mathcal{P}}  \mathbbm{1}(L(M_{n}(P,\textbf{X}_{n}),\hat{M}_{n}(\textbf{X}_{n}))) < 2T_n n^{-\alpha/2})  SP_{\alpha}|\textbf{X}_{n}(dP)  P^{n}_{SP_{\alpha}}(d\textbf{X}_{n}) \\
 & \ \ \ \ \ \ \ \ \ \ \ \ \leq  \limsup_n T_n g(\alpha,s,n)=0. 
\end{align*}

Therefore the first term in \eqref{ciaociao} is equal to zero. Let us consider the second term,
\begin{equation*}
\limsup_n  P^{n}_{SP_{\alpha}} ( K_n(\textbf{X}_{n}) \not \in (n^\alpha/s, s n^\alpha) ).
\end{equation*}
From Theorem 3.8 of Pitman \citet{Pit06}, under the $\alpha$-stable process,  $\frac{K_n(\textbf{X}_{n})}{n^{\alpha}}\rightarrow S_{\alpha}$ almost surely, where $S_{\alpha}$ is a random variable on $\mathbbm{R}_{+}$ distributed according to a Stable distribution of parameter $\alpha$. Therefore, 
\begin{equation} \label{ciaociaociao}
\limsup_n  P^{n}_{SP_{\alpha}} ( K_n(\textbf{X}_{n}) \not \in (n^\alpha/s, s n^\alpha) )= \mathbbm{P} ( S_{\alpha} \not \in (1/s, s) ).
\end{equation}
Taking $s \rightarrow \infty$, \eqref{ciaociaociao} converges to zero, and then so does \eqref{ciaociao}. 

\subsection{Proof of Theorem \ref{thm4}}
In the following, we use the generic notation $C$ to refer to constants that can only depend on $\alpha$ (its value can change from a line to the other).
As in the proof of Theorem \ref{thm3}, let $SP_{\alpha}$ denote the law of a stable process of parameter $\alpha$ and $P^{n}_{SP_{\alpha}}$ the marginal law of the observations under this prior. We can lower bound the minimax risk by the Bayesian risk with prior $SP_{\alpha}$. Indeed,
\begin{align} \label{bau}
  &\inf_{\hat{M}_{n}(\textbf{X}_{n})} \sup_{P\in \mathcal{P}_{PL_\alpha}}\mathbbm{E}_{P}\left( L(\hat{M}_{n}(\textbf{X}_{n}),M_{n}(P,\textbf{X}_{n})) \right) \nonumber \\
& \quad \geq \inf_{\hat{M}_{n}(\textbf{X}_{n})} \int_{\mathcal{P}} \mathbbm{E}_{P}\left( L(\hat{M}_{n}(\textbf{X}_{n}),M_{n}(P,\textbf{X}_{n})) \right) SP_{\alpha}(dP)  \nonumber \\
& \quad = \inf_{\hat{M}_{n}(\textbf{X}_{n})} \int_{[0,1]^{n}} \int_{\mathcal{P}} L(\hat{M}_{n}(\textbf{X}_{n}),M_{n}(P,\textbf{X}_{n}))  SP_{\alpha}|\textbf{X}_{n}(dP) P^{n}_{SP_{\alpha}}(d\textbf{X}_{n}) \nonumber \\
& \quad \geq  \int_{[0,1]^{n}} \inf_{\hat{M}_{n}(\textbf{X}_{n})} \int_{\mathcal{P}} L(\hat{M}_{n}(\textbf{X}_{n}),M_{n}(P,\textbf{X}_{n}))  SP_{\alpha}|\textbf{X}_{n}(dP) P^{n}_{SP_{\alpha}}(d\textbf{X}_{n}) 
\end{align}
From Proposition \ref{Pitman} in Appendix \ref{app}, the posterior distribution of missing mass $M_{n}(P,\textbf{X}_{n})$ under $SP_{\alpha}$ is distributed according to $\text{Beta}(\alpha K_{n}(\textbf{X}_{n}),n-\alpha K_{n}(\textbf{X}_{n}))$. Le $a(\textbf{X}_{n})=\alpha K_{n}(\textbf{X}_{n})$ and $b(\textbf{X}_{n})=n-\alpha K_{n}(\textbf{X}_{n})$, and for easiness of notation we will simply write $a$ and $b$ in the following calculations. The inner integral in \eqref{bau} equals
\begin{align*}
\int_{\mathcal{P}} & L(\hat{M}_{n}(\textbf{X}_{n}),M_{n}(P,\textbf{X}_{n}))  SP_{\alpha}|\textbf{X}_{n}(dP)  \\
 & =\int_{0}^{1}\frac{|\hat{M}_{n}(\textbf{X}_{n})-x|}{x}\frac{x^{a-1}(1-x)^{b-1}}{\mathcal{B}(a,b)}dx \\
& =\frac{\mathcal{B}(a-1,b)}{\mathcal{B}(a,b)}\int_{0}^{1}|\hat{M}_{n}(\textbf{X}_{n})-x|\frac{x^{a-2}(1-x)^{b-1}}{\mathcal{B}(a-1,b)}dx \\
& = \frac{\mathcal{B}(a-1,b)}{\mathcal{B}(a,b)} \mathbbm{E}_{M'}\left(|M'-\hat{M}_{n}(\textbf{X}_{n})|\right)
\end{align*}
where $M'$ is a random variable distributed according to $\text{Be}(a-1,b)$. Plugging this quantity inside \eqref{bau} we find
\begin{align} \label{banana}
&\inf_{\hat{M}_{n}(\textbf{X}_{n})} \sup_{P\in \mathcal{P}_{PL_\alpha}}\mathbbm{E}_{P}\left( L(\hat{M}_{n}(\textbf{X}_{n}),M_{n}(P,\textbf{X}_{n})) \right)  \nonumber \\
&  \qquad\geq\int_{[0,1]^{n}}  \frac{\mathcal{B}(a-1,b)}{\mathcal{B}(a,b)} \inf_{\hat{M}_{n}(\textbf{X}_{n})}  \mathbbm{E}_{M'}\left(|M'-\hat{M}_{n}(\textbf{X}_{n})|\right) P^{n}_{SP_{\alpha}}(d\textbf{X}_{n})  \nonumber \\
& \qquad=  \int_{[0,1]^{n}}  \frac{\mathcal{B}(a-1,b)}{\mathcal{B}(a,b)}  \mathbbm{E}_{M'}\left(|M'-\text{med}(M')|\right) P^{n}_{SP_{\alpha}}(d\textbf{X}_{n}) 
\end{align}
where $\text{med}(M')$ denotes the median of $M'$.
Now, let us denote by $f_{a,b}$ and $m_{a,b}$ the density function and the median of a Beta random variable of parameters $a$ and $b$. We can rewrite the inner expectation in \eqref{banana} as
\begin{align*}
\mathbbm{E}_{M'} & \left(|M' -m_{a-1,b}|\right)  = \int_{0}^{1}|x-m_{a-1,b}|f_{a-1,b}(x)dx \\
& = \int_{0}^{m_{a-1,b}}|x-m_{a-1,b}|f_{a-1,b}(x)dx +  \int_{m_{a-1,b}}^{1}|x-m_{a-1,b}|f_{a-1,b}(x)dx \\
& = - \int_{0}^{m_{a-1,b}}(x-m_{a-1,b})f_{a-1,b}(x)dx +  \int_{m_{a-1,b}}^{1}(x-m_{a-1,b})f_{a-1,b}(x)dx \\
& =   \int_{m_{a-1,b}}^{1}xf_{a-1,b}(x)dx - \frac{1}{2}m_{a-1,b} -  \int_{0}^{m_{a-1,b}}x f_{a-1,b}(x)dx + \frac{1}{2}m_{a-1,b} \\
& =  \int_{m_{a-1,b}}^{1}xf_{a-1,b}(x)dx  -  \int_{0}^{m_{a-1,b}}x f_{a-1,b}(x)dx 
\end{align*}
Therefore, the term inside the integral in \eqref{banana} is
\begin{align*}
&\frac{\mathcal{B}(a-1,b)}{\mathcal{B}(a,b)}  \mathbbm{E}_{M'}\left(|M'-m_{a-1,b}|\right)  \\
&\qquad = \frac{\mathcal{B}(a-1,b)}{\mathcal{B}(a,b)} \int_{m_{a-1,b}}^{1}xf_{a-1,b}(x)dx    - \frac{\mathcal{B}(a-1,b)}{\mathcal{B}(a,b)} \int_{0}^{m_{a-1,b}}x f_{a-1,b}(x)dx  \\
& \qquad= \int_{m_{a-1,b}}^{1}f_{a,b}(x)dx  -  \int_{0}^{m_{a-1,b}} f_{a,b}(x)dx \\
&  \qquad= \int_{m_{a-1,b}}^{m_{a,b}}f_{a,b}(x)dx + \int_{m_{a,b}}^{1}f_{a,b}(x)dx -  \int_{0}^{m_{a,b}} f_{a,b}(x)dx - \int_{m_{a,b}}^{m_{a-1,b}}f_{a,b}(x)dx \\
& \qquad = \int_{m_{a-1,b}}^{m_{a,b}}f_{a,b}(x)dx + \frac{1}{2}  -  \frac{1}{2}- \int_{m_{a,b}}^{m_{a-1,b}}f_{a,b}(x)dx \\
& \qquad = \int_{m_{a-1,b}}^{m_{a,b}}f_{a,b}(x)dx + \int_{m_{a-1,b}}^{m_{a,b}}f_{a,b}(x)dx = 2\int_{m_{a-1,b}}^{m_{a,b}}f_{a,b}(x)dx
\end{align*}
and so we have,
\begin{align} \label{banana2}
\inf_{\hat{M}_{n}(\textbf{X}_{n})} & \sup_{P\in \mathcal{P}_{PL_\alpha}}\mathbbm{E}_{P}\left( L(\hat{M}_{n}(\textbf{X}_{n}),M_{n}(P,\textbf{X}_{n})) \right)  \nonumber \\
&  \geq 2\int_{[0,1]^{n}} \left( \int_{m_{a-1,b}}^{m_{a,b}}f_{a,b}(x)dx \right) P^{n}_{SP_{\alpha}}(d\textbf{X}_{n}) 
\end{align} \vspace{0.4cm}
where we recall that $a=\alpha K_{n}(\textbf{X}_{n})$ and $b=n-\alpha K_{n}(\textbf{X}_{n})$.

From Theorem 3.8 of Pitman \citet{Pit06}, when $P$ is distributed according to the $\alpha$-stable process, $\frac{K_{n}(\textbf{X}_{n})}{n^{\alpha}} \stackrel{d}{\rightarrow}  S_{\alpha}$, where $S_{\alpha}$ is a random variable on $[0,\infty)$ with Stable distribution of parameter $\alpha$. Therefore,  
there exist $ n_{0}$ and two positive bounded values $w_{\alpha}$ and $W_{\alpha}$ such that for all $n>n_{0}$,
 $\mathbbm{P}(\frac{K_{n}(\textbf{X}_{n})}{n^{\alpha}}\in [w_{\alpha},W_{\alpha}])\geq \frac{1}{2}$ and
 $8<2\alpha w_{\alpha}n^{\alpha} \leq 2\alpha W_{\alpha}n^{\alpha} <n-1$.

We are now ready to lower bound \eqref{banana2}. We will make use of some technical lemmas regarding the density and median of the Beta distribution, whose statements and proofs are in Appendix \ref{app}. \eqref{banana2} can be lower bounded by the following quantity
\begin{equation} \label{banana3}
\begin{split}
& 2 \mathbbm{P}(  w_{\alpha}n^{\alpha}  \leq K_{n}(\textbf{X}_{n}) \leq W_{\alpha}n^\alpha) \\
& \qquad\qquad\times \mathbbm{E}_{ P^{n}_{SP_{\alpha}}} \left[ \int_{m_{a-1,b}}^{m_{a,b}}f_{a,b}(x)dx   |  w_{\alpha}n^{\alpha}  \leq K_{n}(\textbf{X}_{n}) \leq W_{\alpha}n^\alpha \right]
 \end{split}
\end{equation}

Given our choice of $ n_{0}$, $w_{\alpha}$ and $W_{\alpha}$,
for $n>n_{0}$, we have that $\mathbbm{P}(  w_{\alpha}n^{\alpha}  \leq K_n \leq W_{\alpha}n^\alpha)\geq \frac{1}{2}$. Recall now that $a=\alpha  K_{n}(\textbf{X}_{n})$ and $b=n- \alpha K_{n}(\textbf{X}_{n})$. Hence, noticing that on the event we are conditioning on, $3<a<b$ and $a < b/2$,  by applying Lemma \ref{lemma2}, we can lower bound \eqref{banana3} by 
\begin{align*}
2 \frac{1}{2}  \mathbbm{E}_{P^{n}_{SP_{\alpha}}} \left[  \frac{C}{\sqrt{a}} |  w_{\alpha}n^{\alpha}  \leq K_{n}(\textbf{X}_{n}) \leq W_{\alpha}n^\alpha  \right],
\end{align*}
for some strictly positive constant $C$. Ultimately this leads to
\begin{align*}
\inf_{\hat{M}_{n}(\textbf{X}_{n})} & \sup_{P\in \mathcal{P}_{RV_\alpha}}\mathbbm{E}_{P}\left( L(\hat{M}_{n}(\textbf{X}_{n}),M_{n}(P,\textbf{X}_{n})) \right)  \\
&\geq C  \mathbbm{E}_{P^{n}_{SP_{\alpha}}} \left[   a^{-1/2} \ | \ w_{\alpha}n^{\alpha}  \leq  K_{n}(\textbf{X}_{n}) \leq W_{\alpha}n^\alpha \right]  \\
&= C \alpha^{-1/2}  \mathbbm{E}_{P^{n}_{SP_{\alpha}}} \left[   K_{n}^{-1/2} \ | \ w_{\alpha}n^{\alpha}   \leq K_{n}(\textbf{X}_{n}) \leq W_{\alpha}n^\alpha \right] \\
& \geq C (w_{\alpha}n^{\alpha})^{-1/2}=Cn^{-\alpha/2}, 
\end{align*}
which provides the lower bound rate for the minimax risk,
\begin{equation*}
\liminf_{n} n^{\alpha/2} \inf_{\hat{M}_{n}(\textbf{X}_{n})}  \sup_{P\in \mathcal{P}_{PL_\alpha}}\mathbbm{E}_{P}\left( L(\hat{M}_{n}(\textbf{X}_{n}),M_{n}(P,\textbf{X}_{n})) \right) >C
\end{equation*}

\subsection{Proof of Proposition \ref{prop2}}

Let $P\in \mathcal{P}_{RV_\alpha}$ and $\ell$ defined as in ($\ref{regular variation}$). In the following, we use the generic notation $C$ and $C'$ to refer to constants that can only depend on $P$ (their values can change from a line to the other).

Here we study the convergence rate under the assumption  of regular variation of the Good-Turing estimator, $ \hat{GT}(\textbf{X}_{n}) = \frac{K_{n,1}(\textbf{X}_{n})}{n}$, proving that
$$ \mathbbm{E}_P \left( L(\hat{GT}(\textbf{X}_{n}), M_{n}(P)) \right) = O(n^{-\alpha/2} \ell(n)^{-1/2}).$$

Let us firs notice that for $a \geq 0,\ b,c > 0$
\begin{equation*}
L(a,c) = \left| \frac{a}{b} \left( \frac{b}{c}-1 \right) + \frac{a}{b}-1 \right| 
\leq L(a,b) + \frac{a}{b} L(b,c).
\end{equation*}
Therefore, we can upper bound the loss of $ \hat{GT}(\textbf{X}_{n}) = \frac{K_{n,1}(\textbf{X}_{n})}{n}$ by
\begin{align*}
&L( \hat{GT}(\textbf{X}_{n}) , M_{n}(P,\textbf{X}_{n}) )
 \\ & \qquad \leq  L (\hat{GT}(\textbf{X}_{n}), \mathbbm{E}_P (\hat{GT}(\textbf{X}_{n}))
+ \frac{\hat{GT}(\textbf{X}_{n})}{\mathbbm{E}_P (\hat{GT}(\textbf{X}_{n}))} L( \mathbbm{E}_P (\hat{GT}(\textbf{X}_{n})), \mathbbm{E}_P (M_{n}(P,\textbf{X}_{n})) \\
&\qquad\qquad\qquad+  \frac{\hat{GT}(\textbf{X}_{n})}{\mathbbm{E}_P (M_{n}(P,\textbf{X}_{n}))} L(\mathbbm{E}_P (M_{n}(P,\textbf{X}_{n})), M_{n}(P,\textbf{X}_{n})),
\end{align*}
and consequently its risk by
\begin{equation}
\begin{split}
& \mathbbm{E}_P (L( \hat{GT}(\textbf{X}_{n}) ,   M_{n}(P,\textbf{X}_{n}) ) ) \leq  \mathbbm{E}_P ( L (\hat{GT}(\textbf{X}_{n}), \mathbbm{E}_P (\hat{GT}(\textbf{X}_{n})))  \\
& \qquad
 + L( \mathbbm{E}_P (\hat{GT}(\textbf{X}_{n})), \mathbbm{E}_P (M_{n}(P,\textbf{X}_{n})) \\
& \qquad\qquad + \mathbbm{E}_P \left( \frac{\hat{GT}(\textbf{X}_{n})}{\mathbbm{E}_P M_{n}(P,\textbf{X}_{n})} L(\mathbbm{E}_P (M_{n}(P,\textbf{X}_{n})), M_{n}(P,\textbf{X}_{n})) \right). \label{upperboundgt}
\end{split}
\end{equation}
We will now separately upper bound the three components of the r.h.s. of this inequality.
Let us first focus on $ L( \mathbbm{E}_P (\hat{GT}(\textbf{X}_{n})), \mathbbm{E}_P (M_{n}(P,\textbf{X}_{n}))$. From Karlin \citet{Kar67} (see also Theorem 4.2 of Ben-Hamou et al. \citet{Ben17}), we know that 
$$ \mathbbm{E}_P (\hat{GT}(\textbf{X}_{n})) \sim \alpha \Gamma(1-\alpha) n^{\alpha-1} \ell(n),$$ 
as $n\rightarrow \infty$ and since  
$0 \leq \mathbbm{E}_P (\hat{GT}(\textbf{X}_{n})) - \mathbbm{E}_P (M_n(P,\textbf{X}_{n})) \leq \frac{1}{n},$
we deduce that
\begin{equation}
L( \mathbbm{E}_P (\hat{GT}(\textbf{X}_{n})), \mathbbm{E}_P (M_{n}(P,\textbf{X}_{n})) \leq C n^{-\alpha} \ell(n)^{-1}\label{ineq_risk_1}.
\end{equation}

Let us now consider the first term in the r.h.s. of \eqref{upperboundgt},
\begin{equation} \label{ineq_risk_02}
\mathbbm{E}_P ( L (\hat{GT}(\textbf{X}_{n}), \mathbbm{E}_P (\hat{GT}(\textbf{X}_{n})))  = \mathbbm{E}_P \left( \left| \frac{K_{n,1}(\textbf{X}_{n})}{\mathbbm{E}_P (K_{n,1}(\textbf{X}_{n}))} - 1\right| \right) .
\end{equation}
As a result of Ben-Hamou et al. \citet{Ben17} Proposition 3.5 (see also the proof of corollary 5.3 of the same paper), for every $\epsilon >0$, we have
$$ P^{n} ( L( K_{n,1}(\textbf{X}_{n}), \mathbbm{E}_P (K_{n,1}(\textbf{X}_{n})) ) \geq \epsilon ) \leq 4 e^{-\epsilon^2 A_n^2}, $$
where 
$$ A_n = \frac{ \mathbbm{E} (K_{n,1}(\textbf{X}_{n})) }{ \sqrt{8 (\mathbbm{E} (K_{n,1}(\textbf{X}_{n})) \vee 2 \mathbbm{E} (K_{n,2}(\textbf{X}_{n}))) + 4/3} }.$$
Hence, we can now bound (\ref{ineq_risk_02}) as follows
\begin{align*}
\mathbbm{E}_P ( L (\hat{GT}(\textbf{X}_{n}),  \mathbbm{E}_P (\hat{GT}(\textbf{X}_{n})) )
&= \mathbbm{E}_P ( L( K_{n,1}(\textbf{X}_{n}), \mathbbm{E}_P (K_{n,1}(\textbf{X}_{n}) ) ) \\
&= \int_{0}^{\infty} P^{n} ( L( K_{n,1}(\textbf{X}_{n}), \mathbbm{E}_P (K_{n,1}(\textbf{X}_{n})) ) \geq \epsilon ) d\epsilon \\
& \leq  4 \int_{0}^{\infty} e^{-\epsilon^2 A_n^2} d\epsilon  =  \frac{4}{A_n} \int_{0}^{\infty} e^{-y^2} dy = C A_n^{-1},
\end{align*}
where we have used the change of variables $y = \epsilon A_n$. Therefore, from the asymptotic behaviors of $\mathbbm{E}_P (K_{n,1}(\textbf{X}_{n}))$ and $\mathbbm{E}_P (K_{n,2}(\textbf{X}_{n}))$ given in Karlin \citet{Kar67} (see also Theorem 4.2 of Ben-Hamou et al. \citet{Ben17}), we conclude that
\begin{equation}
\mathbbm{E}_P ( L (\hat{GT}(\textbf{X}_{n}),  \mathbbm{E}_P (\hat{GT}(\textbf{X}_{n})) ) \leq C n^{-\alpha/2} \ell(n)^{-1/2}\label{ineq_risk_2}.
\end{equation}

Finally, let us look at the third term in \eqref{upperboundgt},
\begin{equation} \label{ineq_risk_03}
\mathbbm{E}_P \left( \frac{\hat{GT}(\textbf{X}_{n})}{\mathbbm{E}_P M_{n}(P,\textbf{X}_{n})} L(\mathbbm{E}_P (M_{n}(P,\textbf{X}_{n})), M_{n}(P,\textbf{X}_{n})) \right).
\end{equation}
Notice that (\ref{ineq_risk_03}) is equal to
\begin{equation} \label{ineq_risk_003}
\mathbbm{E}_P \left(  \frac{\hat{GT}(\textbf{X}_{n})}{M_{n}(P,\textbf{X}_{n})} L(M_{n}(P,\textbf{X}_{n}) , \mathbbm{E}_P (M_{n}(P,\textbf{X}_{n}))) \right),
\end{equation}
and from the Cauchy-Schwarz inequality, we can upper bound \eqref{ineq_risk_003} by
\begin{equation} \label{ineq_risk_0003}
 \sqrt{\mathbbm{E}_P \left( \frac{\hat{GT}(\textbf{X}_{n})^2}{M_{n}(P,\textbf{X}_{n})^2} \right) } \sqrt{ \mathbbm{E}_{P}\left( L(M_{n}(P,\textbf{X}_{n}) , \mathbbm{E}_P (M_{n}(P,\textbf{X}_{n})))^2 \right)}.
\end{equation}
We will first compute the asymptotic behavior of the second term in \eqref{ineq_risk_0003} and then show that the first term is asymptotically bounded. By applying Theorem 3.9 of Ben-Hamou et al. \citet{Ben17} and the asymptotic regimes of Karlin \citet{Kar67}, we obtain that for every $\epsilon>0$, 
$$ P^{n} ( L(M_{n}(P,\textbf{X}_{n}) , \mathbbm{E}_P (M_{n}(P,\textbf{X}_{n})) \geq \epsilon ) \leq 2 e^{-\epsilon^2 B_n}, $$
where 
$$ B_n \leq C n^\alpha \ell(n)$$
(see for example the proof of Corollary 5.3 of Ben-Hamou et al. \citet{Ben17}). Therefore, for all $\epsilon>0$ 
$$P^{n} ( L(M_{n}(P,\textbf{X}_{n}) , \mathbbm{E}_P (M_{n}(P,\textbf{X}_{n})))^2 \geq \epsilon ) \leq 2 e^{-\epsilon B_n} $$
and, following the same reasoning used before, we obtain 
$$ \mathbbm{E}_{P}\left( L(M_{n}(P,\textbf{X}_{n}) , \mathbbm{E}_P (M_{n}(P,\textbf{X}_{n})))^2 \right) \leq C B_n^{-1}$$
which leads to
\begin{equation}
\sqrt{  \mathbbm{E}_{P}\left( L(M_{n}(P,\textbf{X}_{n}) , \mathbbm{E}_P (M_{n}(P,\textbf{X}_{n})))^2 \right) } \leq C  n^{-\alpha/2} \ell(n)^{-1/2}\label{ineq_risk_3}.
\end{equation}

It remains only to prove that the first term in \eqref{ineq_risk_0003} is bounded. First note that
\begin{eqnarray*}
\frac{\hat{GT}(\textbf{X}_{n})}{M_{n}(P,\textbf{X}_{n})} \leq \frac{K_n(\textbf{X}_{n})}{n \sum_{j > K_n(\textbf{X}_{n})} p_{[j]}} 
\end{eqnarray*}
For $t \geq 1$, let us define the function $f$ by $f(t) = \frac{t}{\sum_{j > t} p_{[j]}}$. Noticing that $f(t) \geq 1 > 0$ we can write
\begin{eqnarray}
\frac{\hat{GT}(\textbf{X}_{n})}{M_{n}(P,\textbf{X}_{n})} \leq \frac{f(K_n(\textbf{X}_{n}))}{f(\mathbbm{E} (K_n(\textbf{X}_{n})))} \frac{f(\mathbbm{E} (K_n(\textbf{X}_{n})))}{n}\label{ineq_risk_032}.
\end{eqnarray}
Denoting by $\ell^{\frac{1}{\alpha} \#}$ the de Bruijn conjugate of $\ell^{\frac{1}{\alpha}}$ (see subsection 1.5.7 of Bingham et al. \citet{Bin87} for a definition), Proposition 23 of Gnedin et al. \citet{Gne07} implies that
$$ f(t) \sim C t^{\frac{1}{\alpha}} \ell^{\frac{1}{\alpha} \#}(t^{\frac{1}{\alpha}}), $$
which in turns implies that $f^2$ is regularly varying with index $2/\alpha$.  Since $f$ is non decreasing, it is bounded on any set of the form $[1,T]$. Therefore, we can apply Potter's Theorem (Theorem 1.5.6, Bingham et al. \citet{Bin87}) 
to obtain
$$ \frac{f(K_n(\textbf{X}_{n}))^2}{f(\mathbbm{E}_P (K_n(\textbf{X}_{n})))^2} \leq C \left[\left(\frac{K_n(\textbf{X}_{n})}{\mathbbm{E}_P (K_n(\textbf{X}_{n}))}\right)^{\frac{1}{\alpha}} + \left(\frac{K_n(\textbf{X}_{n})}{\mathbbm{E}_P (K_n(\textbf{X}_{n}))}\right)^{\frac{3}{\alpha}}\right] + C'.$$
Following the same reasoning used before,  we can show that for all $\eta > 1$
\begin{equation}\label{ineq_K_3}
 \lim_{n\rightarrow +\infty} \mathbbm{E}_P [L(K_n(\textbf{X}_{n}), \mathbbm{E}_P (K_n(\textbf{X}_{n})))^\eta] = 0,
\end{equation}
and, thanks to the elementary inequality $x/|x-1| \leq 2$ for $x \geq 2$, it follows that, for all $\eta > 1$,
$$  \left(\frac{K_n(\textbf{X}_{n})}{\mathbbm{E}_P (K_n(\textbf{X}_{n}))}\right)^\eta \leq 2^\eta + 2^\eta L(K_n(\textbf{X}_{n}), \mathbbm{E}_P (K_n(\textbf{X}_{n})))^\eta.$$
As a consequence of this last inequality, along with (\ref{ineq_K_3}), for all $\eta>1$ we obtain
$$  \mathbbm{E}_P \left( \left(\frac{K_n(\textbf{X}_{n})}{\mathbbm{E}_P (K_n(\textbf{X}_{n}))}\right)^\eta  \right) \leq C $$
from which we get that
\begin{equation}\label{ineq_K_31}
\mathbbm{E}_P \left( \left(\frac{f( K_n(\textbf{X}_{n}))}{f( \mathbbm{E}_P (K_n(\textbf{X}_{n})))}\right)^\eta  \right) \leq C.
\end{equation}

Besides, since $ \mathbbm{E}_P (K_n(\textbf{X}_{n}))^{\frac{1}{\alpha}} \sim n \ell^{\frac{1}{\alpha}}(n)$, which diverges to infinity, the uniform convergence theorem for slowly varying functions (Theorem 1.2.1, Bingham et al. \citet{Bin87}) gives that
$$ \ell^{\frac{1}{\alpha}\#}(\mathbbm{E}_P (K_n(\textbf{X}_{n}))^{\frac{1}{\alpha}} ) \sim \ell^{\frac{1}{\alpha}\#} (n \ell^{\frac{1}{\alpha}}(n)).$$
As a consequence of this and of the asymptotic properties of $f$, we obtain that
$$ \frac{f(\mathbbm{E}_P (K_n(\textbf{X}_{n})))}{n} \sim \ell^{\frac{1}{\alpha}}(n) \ell^{\frac{1}{\alpha}\#} (n \ell^{\frac{1}{\alpha}}(n)), $$
which, from the definition of the de Bruijn conjugate, in turn gives 
$$ \frac{f(\mathbbm{E}_P (K_n(\textbf{X}_{n})))}{n} \sim 1, $$
  and then
\begin{equation}
\frac{f(\mathbbm{E}_P (K_n(\textbf{X}_{n})))^2}{n^2} \leq C \label{ineq_K_32}.
\end{equation}
From (\ref{ineq_risk_032}), (\ref{ineq_K_31}) and (\ref{ineq_K_32}) together, we finally obtain 
$$  \mathbbm{E}_P \left( \frac{\hat{GT}(\textbf{X}_{n})^2}{M_{n}(P,\textbf{X}_{n})^2} \right) \leq C,$$
which together with (\ref{ineq_risk_1}), (\ref{ineq_risk_2}) and (\ref{ineq_risk_3}) concludes the proof.

%%%%%%%%%%%%%%%%%%%%%%%%%%%%%%%%
%%%%%%%%%%%%%%%%%%%%%%%%%%%%%%%%
%%%%%%%%%%%%%%%%%%%%%%%%%%%%%%%%
%%%%%%%%%%%%%%%%%%%%%%%%%%%%%%%%

%\appendix
%\normalsize
\section{Appendix}\label{app}

\begin{proposition} \label{prop:measurable_missing}
The missing mass $M_{n}(P,\textbf{X}_{n})$ is a jointly measurable mapping.
\end{proposition}
\textsc{proof}
Recall that $\mathcal{P}$ is endowed with the smallest $\sigma$-algebra making the mappings $P\mapsto P(A)$ measurable for every $A\in \mathcal{B}([0,1])$. This is also the Borel $\sigma$-algebra generated by the weak convergence topology, which can be induced by the \textit{bounded Lipschitz metric} (see Appendix A of Ghosal and Van der Vaart \citet{Gho17}), defined as
$$ d_{BL}(P,Q) = \sup_{ \|f\|_{\mathcal{C}^1} \leq 1} \left|\int f dP - \int f dQ \right|$$
where the supremum is over all real functions satisfying $ |f(x)-f(y)| \leq |x-y| $ for any $x,y \in [0,1]$. $[0,1]^n$ is endowed with the Euclidean topology, which can be induced by the $\ell_\infty$ norm.  \par

Let us consider $O_{n}(P,\textbf{X}_{n}) = 1-M_{n}(P,\textbf{X}_{n})$ and let us define for any $\eta > 0$, 

\begin{equation*}
 f_{\eta,X_i}:x \mapsto \max(0, \eta -|X_i-x|).
\end{equation*}

Also let $f_{\eta,\textbf{X}_{n}} = \max_i f_{\eta,X_i}$, which is  $1-$Lipschitz function, since all $f_{\eta,X_i}$ are $1-$Lipschitz. Now, let $O_{\eta,n}$ be defined as follows

\begin{equation*}
 O_{\eta,n}(P,\textbf{X}_{n}) = \frac{1}{\eta}\sum_{j\geq1} p_j f_{\eta,\textbf{X}_{n}}(\theta_j) = \int \frac{1}{\eta}f_{\eta,\textbf{X}_{n}} dP.
\end{equation*}

Let $(P,\textbf{X}_{n}) \in  \mathcal{P}\times [0,1]^n$. We have that
\begin{equation} \label{missO}
\lim_{\eta \rightarrow 0} O_{\eta,n}(P,U_{1:n}) = O_{n}(P,U_{1:n}).
\end{equation} 
Indeed, notice that for any $x \in [0,1]$, 

\begin{equation*}
\frac{f_{\eta,\textbf{X}_{n}}(x)}{\eta} > 0 \Leftrightarrow \exists i \in \{1,..,n\},\ |X_i - x | <\eta. 
\end{equation*}

Hence, if $x \not \in \textbf{X}_{n}, \lim_{\eta \rightarrow 0} \frac{f_{\eta,\textbf{X}_{n}}(x)}{\eta} = 0 .$ Besides, $\frac{f_{\eta,\textbf{X}_{n}}(X_i)}{\eta} = 1 $ for any $\eta$ and $i$. Finally, since $\frac{f_{\eta,\textbf{X}_{n}}(X_i)}{\eta} \leq 1$, dominated convergence theorem gives \eqref{missO}.\\

Let $\epsilon > 0$, and take $\textbf{X}_{n},\textbf{Y}_{n} \in [0,1]^{n}$ such that $\| \textbf{X}_{n}-\textbf{Y}_{n}\|_\infty \leq \eta \epsilon/2$. Take $P,Q \in \mathcal{P}$ such that $d_{BL}(P,Q) \leq \eta \epsilon/2$. Now,  for any $x \in [0,1]$, 
\begin{equation*} 
|f_{\eta,\textbf{X}_{n}}(x) - f_{\eta,\textbf{Y}_{n}}(x)| \leq \eta \epsilon/2.
\end{equation*}
Indeed, suppose for instance $f_{\eta,\textbf{X}_{n}}(x) \geq f_{\eta,\textbf{Y}_{n}}(x)$. Suppose $f_{\eta,\textbf{X}_{n}}(x) > 0$. Now consider $X_i$ the closest point to $x$. We have that 
\begin{equation*}
f_{\eta,\textbf{X}_{n}}(x) = \eta - | X_i - x |  \leq \eta - | Y_i - x | + \eta\epsilon/2 = f_{\eta,\textbf{Y}_{n}}(x)+ \eta\epsilon/2.
\end{equation*}

Finally, let us compute the distance between the two images,
\begin{align*}
\eta & \left| O_{\eta,n}(P,\textbf{X}_{n})  - O_{\eta,n}(Q,\textbf{Y}_{n}) \right| =  \left|\int f_{\eta,\textbf{X}_{n}} dP - \int f_{\eta,\textbf{Y}_{n}} dQ \right| \\
& \ \ \ \ \ \ \ \ \leq  \left|\int f_{\eta,\textbf{X}_{n}} dP - \int f_{\eta,\textbf{X}_{n}} dQ \right|  + \left|\int f_{\eta,\textbf{X}_{n}} dQ - \int f_{\eta,\textbf{Y}_{n}} dQ \right| \\
&  \ \ \ \ \ \ \ \ \leq \left|\int f_{\eta,\textbf{X}_{n}} dP - \int f_{\eta,\textbf{X}_{n}} dQ \right| + \| f_{\eta,\textbf{X}_{n}} - f_{\eta,\textbf{Y}_{n}} \|_\infty \leq \eta\epsilon
\end{align*}
which gives
$$ | O_{\eta,n}(P,\textbf{X}_{n}) - O_{\eta,n}(Q,\textbf{Y}_{n})| \leq \epsilon. $$
Therefore $O_{\eta,n}$ is continuous and hence measurable. And finally $M_n$ is measurable as limit and sum of measurable functions.
\qed

\begin{proposition} \label{propPitman}
Let $\textbf{X}_{n}=(X_{1},\ldots,X_{n})$ be a sample such that $X_{i}|P\overset{iid}{\sim}P$ for all $1\leq i\leq n$. Then,
\begin{enumerate}
\item if  $P \sim \text{DP}(\gamma)$, where $\gamma(d\theta)=\mathbbm{1}(0<d\theta<1)$, then $M_{n}(\textbf{X}_{n},P)|\textbf{X}_{n} \sim \text{Beta}(1,n)$
\item if $P \sim \text{SP}_{\alpha}$,  then $M_{n}(\textbf{X}_{n},P)|\textbf{X}_{n} \sim \text{Beta}(\alpha K_{n}(\textbf{X}_{n}),n - \alpha K_{n}(\textbf{X}_{n}))$
\end{enumerate}

\end{proposition}
\textsc{proof}
We are going the derive the posterior distribution of  $M_{n}(\textbf{X}_{n},P)$ when $P$ is distributed according to a Pitman-Yor process (Perman et al. \citet{Per92} and Pitman and Yor \citet{Pit97}), $P \sim \text{PY}(\eta,\alpha)$, with $\alpha < 1$ and $\eta > -\alpha$. Point i) in the statement is the particular case  $\text{PY}(1,0)$, while point ii) corresponds to $\text{PY}(0,\alpha)$.

From Corollary 20 of Pitman \citet{Pit96}, the posterior distribution of $P$ given $\textbf{X}_{n}$ under a Pitman-Yor process satisfies the following distributional equality,
\begin{equation} \label{Pitman}
P|\mathbf{X}_{n}\overset{d}{=}\sum_{i=1}^{K_{n}}w_{i}\delta_{X_{i}^{*}}+w_{0}\tilde{P}
\end{equation}
where $(X_{1}^{*},\ldots,X_{K_{n}}^{*})$ are $K_{n}$ the distinct values in the sample $\mathbf{X}_{n}$ and having multiplicities $(n_{1},\ldots,n_{K_{n}})$, $w=\left(w_{0},w_{1},\ldots,w_{K_{n}}\right)$ is a random vector distributed according to a Dirichlet distribution
 $\text{Dir}\left(\eta+K_{n}\alpha,n_{1}-\alpha,\ldots,n_{K_{n}}-\alpha\right)$ and $\tilde{P}\sim \text{PY}(\alpha,\eta+K_{n}\alpha)$ independent of $w$. 
Therefore, 
\begin{equation} \label{Pitman1}
M_{n}(\textbf{X}_{n},P)=P(\{\mathbf{X}_{n}\}^{c})|\mathbf{X}_{n}\overset{d}{=}\sum_{i=1}^{K_{n}}w_{i}\delta_{X_{i}^{*}}(\{\mathbf{X}_{n}\}^{c})+w_{0}\tilde{P}(\{\mathbf{X}_{n}\}^{c})
\end{equation}
The point masses in \eqref{Pitman1} are all equal to zero, while $\tilde{P}(\{\mathbf{X}_{n}\}^{c})=1$ since the base measure of $\tilde{P}$ is diffuse. Therefore, $M_{n}(\textbf{X}_{n},P)\overset{d}{=}w_{0}$ and $w_{0}$ is distributed according to $\text{Beta}(\eta+K_{n}\alpha,n-\alpha K_{n})$ from the aggregation property of the Dirichlet distribution.
\qed

\begin{lemma} \label{lemma1}
Let $f_{a,b}$ denote the density of the $\text{Beta}(a,b)$ distribution. Then, there exists $n_{0}\in \mathbbm{N}$ such that for all $b>a>n_{0}$, we have
\begin{equation} \label{beta11}
\sup_{x\in [0,1]} f_{a,b}(x) < 8 (a+b)^{3/2} a^{-1/2} b^{-1/2}
\end{equation}
\end{lemma}
\textsc{proof} 
Take first $a,b > 2$. The mode of the beta distribution is $\frac{a-1}{a+b-2} $, therefore
\begin{align}
\sup_{x\in [0,1]} f_{a,b}(x) & = f_{a,b} \left(\frac{a-1}{a+b-2} \right)  \leq \frac{1}{\mathcal{B}(a,b)} \frac{a^{a-1} b^{b-1}}{ (a+b-2)^{a+b-2} } \label{bananas2}
\end{align}
From the Stirling's formula, there exists $n_0$ such that for $a,b > n_0$,
$$ \frac{1}{\mathcal{B}(a,b)} \leq \frac{ (a+b)^{a+b-1/2} }{ a^{a-1/2} b^{b-1/2}}.$$
Plugging this quantity into \eqref{bananas2}, we have
$$ \sup_{x\in [0,1]} f_{a,b}(x) \leq a^{-1/2} b^{-1/2} (a+b)^{3/2} \frac{1}{(1-\frac{2}{a+b})^{a+b}}.  $$
Finally, since for $n_0$ large enough $\frac{1}{(1-\frac{2}{a+b})^{a+b}} \leq 8$, we have \eqref{beta11}.
\qed

\begin{lemma}  \label{lemma2}
There exists a constant $C$ such that for any $a,b > 3$ such that $a<b/2$, 
\begin{equation}
\int_{m_{a-1,b}}^{m_{a,b}}f_{a,b}(x)dx \geq \frac{C}{\sqrt{a}}
\end{equation}
\end{lemma}
\textsc{proof}
In the following, we use the generic notation $C$ to refer to universal constants. The value of $C$ can change from a line to another. To simplify the notations, let $n = a+b$ and denote $I_x(a,b)$ the normalized incomplete Beta function defined as
$$ I_x(a,b) =  \int_{0}^{x}f_{a,b}(x)dx.$$
It is well-known that 
$$ I_x(a+1,b) = I_x(a,b) - \frac{x^a (1-x)^b}{a \mathcal{B}(a,b)}, $$
(such a result can be obtain using integration by part). Successively applying this result for $x = m_{a-1,b}$ and using the definition of the median we deduce that
\begin{eqnarray*}
I_{m_{a-1,b}}(a,b) &=& I_{m_{a-1,b}}(a-1,b) - \frac{m_{a-1,b}^{a-1} (1-m_{a-1,b})^b}{(a-1) \mathcal{B}(a-1,b)} \\
&=& 1/2 - \frac{m_{a-1,b}^{a-1} (1-m_{a-1,b})^b}{(a-1) \mathcal{B}(a-1,b)} \\
&=& I_{m_{a,b}}(a,b) - \frac{m_{a-1,b}^{a-1} (1-m_{a-1,b})^b}{(a-1) \mathcal{B}(a-1,b)},
\end{eqnarray*}
which in turn leads to
\begin{equation}
\int_{m_{a-1,b}}^{m_{a,b}}f_{a,b}(x)dx =  I_{m_{a,b}}(a,b) - I_{m_{a-1,b}}(a,b) = \frac{m_{a-1,b}^{a-1} (1-m_{a-1,b})^b}{(a-1) \mathcal{B}(a-1,b)}.\label{ineq_diff_b}
\end{equation}

From the Stirling formula applied to the Beta function, we know that for all $a,b>1$,
\begin{equation*}
\mathcal{B}(a,b) \leq C \frac{a^{a-1/2} b^{b-1/2}}{(a+b)^{a+b-1/2}}=C\frac{a^{a-1/2} b^{b-1/2}}{n^{n-1/2}}, 
\end{equation*}
which leads to 
\begin{equation}
\int_{m_{a-1,b}}^{m_{a,b}}f_{a,b}(x)dx \geq C \frac{(n-1)^{n-3/2}}{(a-1)^{a-1/2} b^{b-1/2}}m_{a-1,b}^{a-1} (1-m_{a-1,b})^b\label{ineq_diff_b2}
\end{equation}
when plugging in (\ref{ineq_diff_b}).

Now, since $a-1 < b$, the mode-median-mean inequality (see Payton et al. \citet{Pay89}) gives that
$$ \frac{a-2}{n-3} \leq m_{a-1,b} \leq \frac{a-1}{n-1},$$
from which we deduce
$$ m_{a-1,b}^{a-1} \geq \frac{(a-2)^{a-1}}{(n-3)^{a-1}} \geq \frac{(a-2)^{a-1}}{(n-1)^{a-1}}$$
and
\begin{equation*}
(1-m_{a-1,b})^b \geq \frac{b^b}{(n-1)^b}.
\end{equation*}
Now, together with (\ref{ineq_diff_b2}), the previous two inequalities yield
\begin{equation*}
\int_{m_{a-1,b}}^{m_{a,b}}f_{a,b}(x)dx \geq C \frac{(n-1)^{n-3/2}}{(a-1)^{a-1/2} b^{b-1/2}}\frac{(a-2)^{a-1} b^b }{(n-1)^{n-1}}  \geq C \sqrt{\frac{b}{a n}}
\end{equation*}
where we used the fact that for $x > 2$, $x^x \geq (x-1)^x \geq C x^x$. We then conclude by noticing that $a < b/2$ implies that $b/n < 2/3$.
\qed

%%%%%%%%%%%%%%%%%%%%%%%%%%%%%%%%
%%%%%%%%%%%%%%%%%%%%%%%%%%%%%%%%
%%%%%%%%%%%%%%%%%%%%%%%%%%%%%%%%
%%%%%%%%%%%%%%%%%%%%%%%%%%%%%%%%

\end{document}